\documentclass[12pt]{article}
\usepackage{amssymb,latexsym,amsmath}
\usepackage{graphics}                 
\usepackage{color}                    
\usepackage{hyperref}                 
\usepackage[all]{xy}

\begin{document}
\pagestyle{plain} \headheight=5mm \topmargin=-5mm 

\title{The Lawson homology and Deligne-Beilinson cohomology for Fulton-MacPherson configuration spaces }
\author{ Wenchuan Hu and Li Li
}

\maketitle
\newtheorem{definition}{Definition}[section]
\newtheorem{theorem}{Theorem}[section]
\newtheorem{proposition}{Proposition}[section]
\newtheorem{lemma}{Lemma}[section]
\newtheorem{corollary}{Corollary}[section]
\newtheorem{question}{Question}
\newtheorem{remark}{Remark}[section]

\def\s{\section}
\def\ss{\subsection}

\def\nn{\nonumber}
\def\bp{{\bf Proof.}\hspace{2mm}}
\def\qe{\hfill$\Box$}
\def\lj{\langle}
\def\rj{\rangle}
\def\ox{\mbox{}}
\def\lb{\label}
\def\mod{\;{\rm mod}\;}
\def\exp{{\rm exp}\;}
\def\Lie{{\rm Lie}}
\def\dim{{\rm dim}}
\def\im{{\rm im}\;}
\def\Supp{{\rm Supp}\;}
\def\Sp{{\rm Sp}\;}
\def\ind{{\rm ind}\;}
\def\rank{{\rm rank}\;}
\def\deg{{\rm deg}}
\def\ker{{\rm ker}\;}

\def\H{{\bf H}}
\def\K{{\rm K}}
\def\R{{\bf R}}
\def\C{{\bf C}}
\def\Z{{\mathbb{Z}}}
\def\ch{{\rm ch}}
\def\P{{\rm P}}

\begin{abstract}
In this paper, we compute the Lawson homology groups and Deligne-Beilinson
cohomology groups for Fulton-MacPherson configuration spaces. The
explicit formulas are given.
\end{abstract}

\begin{center}{\bf \tableofcontents}\end {center}

\section{Introduction}

\hskip .2in In this paper, all varieties are defined over
$\mathbb{C}$. Let $X$ be an $d$-dimensional projective variety. Let
${\cal Z}_p(X)$ be the space of algebraic $p$-cycles on $X$.

The \textbf{Chow group} ${\rm Ch}_p(X)$  of p-cycles is defined by
${\cal Z}_p(X)$ modulo the rational equivalence. For general
background on Chow groups, the reader is referred to Fulton's  book
\cite{Fulton}.

The \textbf{Lawson homology} $L_pH_k(X)$ of $p$-cycles is defined by
$$L_pH_k(X) := \pi_{k-2p}({\cal Z}_p(X)) \quad {\rm for}\quad k\geq 2p\geq 0,$$
where ${\cal Z}_p(X)$ is provided with a natural topology (cf.
\cite{Friedlander1}, \cite{Lawson1}). For general background on
Lawson homology, the reader is referred to \cite{Lawson2}.

It is convenient to extend the definition of Lawson homology by
setting
$$ L_pH_k(X)=L_0H_k(X), \quad {\rm if} \quad p<0.$$

It was proved in \cite{author1} that, for any smooth projective
variety $X$, the formula on  Lawson homology for a blowup holds:

\begin{theorem}[\cite{author1}]
 { Let $X$ be smooth projective
variety and $Y\subset X$ be a smooth subvariety of codimension
$r\geq 2$. Let $\sigma:\tilde{X}_Y\rightarrow X$ be the blowup of
$X$ along $Y$, $\pi:D:=\sigma^{-1}(Y)\rightarrow Y$ the natural map,
and $i:D\rightarrow \tilde{X}_Y$ the exceptional divisor of the
blowup. Then for integers $p$, $k$ with $k\geq 2p\geq 0$, there is
an isomorphism
  $$
\begin{array}{cc}
  &I_{p,k}: \bigg\{\bigoplus_{1\leq j \leq
r-1}L_{p-j}H_{k-2j}(Y) \bigg\}\oplus L_pH_k(X)\cong L_pH_k(\tilde{X}_Y).
\end{array}
  $$
 }
\end{theorem}

 Now we give minimal notations for the Fulton-MacPherson
configuration spaces enough for stating the main theorem (see
section \ref{sec:2.2} for a construction of the Fulton-MacPherson
configuration spaces by a sequence of blowups.)

Let $X$ be a smooth projective variety of dimension $d$ and let
$n\ge 1$ be an integer. Consider the cartesian product $X^n:=X\times\cdots \times X$
of $n$ copies of $X$. Denote
by $\Delta_I$ the diagonal in $X^n$ where $x_i=x_j$ if $i,j\in I$.

The \emph{configuration space} $F(X,n)$ is the complement of all
diagonals in $X^n$, i.e.,

$$F(X,n)=X^n\setminus \cup_{|I|\ge 2}\Delta_I=\{(x_1,\dots,x_n)\in X^n: x_i\neq x_j, \forall i\neq j\}.$$

For each subset $I\in [n]:=\{1,\dots, n\}$ with at least two
elements, denote by $\textrm{Bl}_\Delta(X^I)$ the blowup of the
corresponding cartesian product $X^I$ along its small diagonal. In
\cite{FultonM}, Fulton and MacPherson have given the definition of
their compactification $X[n]$ as follows.

\begin{theorem}[Fulton-MacPherson]
{ The closure of the natural locally closed
embedding
$$i:F(X,n)\hookrightarrow X^n\times \prod_{|I|\ge 2}\textrm{Bl}_\Delta(X^I)$$
is smooth, and the boundary is a simple normal crossing divisor. The
closure is called the Fulton-MacPherson configuration space, denoted
by $X[n]$.}
\end{theorem}

 We call two subsets $I, J\subseteq[n]:=\{1,2,\dots,n\}$ are \emph{overlapped} if $I\cap J$ is
a nonempty proper subset of $I$ and of $J$.

A \emph{nest} $\mathcal{S}$ is a set of subsets of $[n]$ such that
any two elements $I\neq J\in \mathcal{S}$ are not overlapped, and
all singletons $\{1\},\dots,\{n\}$ are in $\mathcal{S}$. Notice that
the nest defined here, unlike the one defined in \cite{FultonM},
contains singletons.

Given a nest $\mathcal{S}$, define
$\mathcal{S}^\circ=\mathcal{S}\setminus
\big{\{}\{1\},\dots,\{n\}\big{\}}$. In the description of nests by
forests below, $\mathcal{S}^\circ$ corresponds to the forest
$\mathcal{S}$ cutting of all leaves.

A nest $\mathcal{S}$ naturally corresponds to a not necessarily
connected tree (which is also called a \emph{forest} or a
\emph{grove}), each node of which is labeled by an element in
$\mathcal{S}$. For example, the following forest corresponds to a
nest $\mathcal{S}=\{1,2,3,23,123\}$.

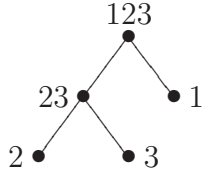
\begin{figure}[h]
\hskip50mm\xy
    (150,-20)="m";
    "m"+(12,16)="123"  *{\bullet};
    "m"+(6,8)="23"  *{\bullet};
    "m"+(18,8)="1"  *{\bullet};
    "m"+(0,0)="2"  *{\bullet};
    "m"+(12,0)="3"  *{\bullet};
    "m"+(12,19)  *{123};
    "m"+(2,8)  *{23};
    "m"+(21,8)  *{1};
    "m"+(-3,0)  *{2};
    "m"+(15,0)  *{3};
    "123";         "23"**\dir{-};
    "123";         "1"**\dir{-};
    "23";         "2"**\dir{-};
    "23";         "3"**\dir{-};
\endxy
\end{figure}

Denote by $c(\mathcal{S})$ the number of connected components of the
forest, i.e., the number of maximal elements of $\mathcal{S}$.
Denote by $c_I(\mathcal{S})$ (or $c_I$ if no ambiguity arise) the
number of maximal elements of the set $\{J\in\mathcal{S}|J\subsetneq
I\}$, i.e. the number of sons of the node $I$. In the above example,
$c(\mathcal{S})=1$, $c_{123}=c_{23}=2$.

For a nest $\mathcal{S}\neq\{\{1\},\dots,\{n\}\}$ (i.e.
$\mathcal{S}^\circ\neq\emptyset$), define a set $M_\mathcal{S}$ of
lattice points in the integer lattice
$\mathbb{Z}^{\mathcal{S}^\circ}$ as follows
$$M_\mathcal{S}:=\big{\{}\underline{\mu}=\{\mu_I\}_{I\in\mathcal{S}^\circ}:
1\le \mu_I\le d(c_I-1)-1 \big{\}}.$$ (Recall that $d=\dim X$,
$c_I=c_I(\mathcal{S})$) and define
$||\underline{\mu}||:=\sum_{I\in\mathcal{S}^\circ} \mu_I$, $\forall
\underline{\mu}\in M_\mathcal{S}$.

For $\mathcal{S}=\{\{1\},\dots,\{n\}\}$, assume
$M_\mathcal{S}=\{\underline{\mu}\}$ with $\|\underline{\mu}\|=0$.

It was proved in \cite{author2} that, for any smooth projective
variety $X$ the following holds:

\begin{theorem}[\cite{author2}]\label{li thm}
{Let $X$ be a smooth projective variety defined
over $\mathbb{C}$. Then for each $p\ge 0$, there is an isomorphism
of Chow groups:
$${\rm Ch}_p(X[n])\cong\bigoplus_{\mathcal{S}}\bigoplus_{\underline{\mu}\in
M_{\mathcal{S}}} Ch_{p-||\underline{\mu}||}(X^{c(\mathcal{S})}) .$$
where $\mathcal{S}$ runs through all nests of $[n]$.}
\end{theorem}

The first main result in this paper is the following
\begin{theorem}
{\label{sec:M1} Let $X$ be a smooth projective variety defined
over $\mathbb{C}$. Then for each pair of integers $p,k$, $k\geq
2p\geq 0$, there is an isomorphism of Lawson homology groups:
$$L_pH_k(X[n])\cong\bigoplus_{\mathcal{S}}\bigoplus_{\underline{\mu}\in
M_{\mathcal{S}}}
L_{p-||\underline{\mu}||}H_{k-2||\underline{\mu}||}(X^{c(\mathcal{S})})
.$$ where $\mathcal{S}$ runs through all nests of $[n]$.
}
\end{theorem}

\begin{remark}
When $p=0$, Theorem \ref{sec:M1} reduces to the formula of singular 
homology groups with integer coefficient for $X[n]$. 
In particular, the integer singular homology of $X[n]$ depends only
the integer singular homology of $X$.
\end{remark}

As a corollary, we have the following more explicit formula:

\begin{corollary}
{Let $X$ be a smooth projective variety defined over
$\mathbb{C}$. Then for each pair of integers $p,k$, $k\geq 2p\geq
0$, there is an isomorphism of Lawson homology groups:
$$
L_pH_k(X[n])\cong\bigoplus_{\substack{1\le m\le n\\
0\leq i\leq p}} L_{p-i}H_{k-2i}(X^m)^{\oplus
[\frac{x^it^n}{n!}]\frac{N^m}{m!}}.
$$ where $N$  and $\oplus
[\frac{x^it^n}{n!}]\frac{N^m}{m!}$ are defined in (\ref{eq:frac}).
}
\end{corollary}

Let $X$ be a complex manifold of complex dimension $d$. Let
$\Omega_X^k$ the sheaf of holomorphic $k$-form on $X$. The
\textbf{Deligne complex of level p} is the complex of sheaves
$$\underline{\Z}_{\cal{D}}(p):0\rightarrow\Z\stackrel{(2i\pi)^p}{\rightarrow}\Omega_X^0\rightarrow
\Omega_X^1\rightarrow\Omega_X^2\rightarrow\cdots\rightarrow\Omega_X^{p-1}\rightarrow
0
$$
The \textbf{Deligne-Beilinson cohomology} of $X$ in level $p$ we
mean the hypercohomology of this complex:
$$ H^*_{\cal{D}}(X,\Z(p)):= \mathbb{H}^*(X,\underline{\Z}_{\cal{D}}(p))
$$

For Deligne-Beilinson cohomology $H^k_{\cal{D}}(-,\Z(p))$, we obtain the
following result:

\begin{theorem}
{Let $X$ be a smooth projective variety defined over $\mathbb{C}$.
Then for each pair of integers $p$, $k$,  there is an isomorphism of
Deligne-Beilinson cohomology groups:
$$
H^k_{\mathcal{D}}(X[n],\Z(p))\cong\bigoplus_{\mathcal{S}}\bigoplus_{\underline{\mu}\in
M_{\mathcal{S}}}
H_{\mathcal{D}}^{k-2||\underline{\mu}||}(X^{c(\mathcal{S})},\Z({p-||\underline{\mu}||})).
$$
}
\end{theorem}

\medskip
The main tools used to prove the main result are: The formula on the
Lawson homology for a blowup proved  in \cite{author1} and the
method in computing the Chow groups of the Fulton-MacPherson
configuration space in \cite{author2}.

 \section{Some fundamental
materials}

\ss{Lawson homology} {\hskip .2 in} Recall that for a morphism
$f:U\rightarrow V$ between  projective varieties, there exist
induced homomorphism $f_*: L_pH_k(U)\rightarrow L_pH_k(V)$ for all
$k\geq 2p\geq 0$. Furthermore, it has been shown by C. Peters
\cite{Peters} that if $U$ and $V$ are smooth and projective, there
are Gysin ``wrong way" homomorphism $f^*: L_pH_k(V)\rightarrow
L_{p-c}H_{k-2c}(U)$, where $c=\dim(V)-\dim(U)$.

\medskip
Let $X$ be a smooth projective variety and $i_0:Y\hookrightarrow X$
a smooth subvariety of codimension $r$. Let
$\sigma:\tilde{X}_Y\rightarrow X$ be the blowup of $X$ along $Y$,
$\pi:D=\sigma^{-1}(Y)\rightarrow Y$ the natural map, and
$i:D=\sigma^{-1}(Y)\hookrightarrow \tilde{X}_Y$ the exceptional
divisor of the blowing up. Set $U\equiv X-Y\cong \tilde{X}_Y - D$.
Denote by $j_0$ the inclusion $U\subset X$ and $j$ the inclusion
$U\subset \tilde{X}_Y$. Note that $\pi:D=\sigma^{-1}(Y)\rightarrow
Y$ makes $D$ into a projective bundle of rank $r-1$, given precisely
by $D=\P(N_{Y/X})$ and we have (cf. \cite{Voisin2}, pg.271)

$$ {\mathcal{O}}_{\tilde{X}_Y}(D)|_{D}=
{\mathcal{O}}_{\P(N_{Y/X})}(-1).
$$

Denote by $h$ the class of ${\mathcal{O}}_{\P(N_{Y/X})}(-1)$ in
${\rm Pic}(D)$. We have $h=-D_{|D}$ and
$-h=i^*i_*:L_qH_m(D)\rightarrow L_{q-1}H_{m-2}(D)$ for $0\leq 2q\leq
m$ ([\cite{FG}, Theorem 2.4], [\cite{Peters}, Lemma 11]). The last
equality can be equivalently regarded as a Lefschetz operator
\begin{eqnarray*}
-h=i^*i_*:L_qH_m(D)\rightarrow L_{q-1}H_{m-2}(D), \quad 0\leq 2q\leq
m.
\end{eqnarray*}

\medskip
The proof of the main result are based on the following Theorems:

\begin{theorem}[Lawson homology for a blowup]
{\label{sec:bulh} Let $X$ be
smooth projective manifold and $Y\subset X$ be a smooth subvariety
of codimension r. Let $\sigma:\tilde{X}_Y\rightarrow X$ be the
blowup of $X$ along $Y$, $\pi:D=\sigma^{-1}(Y)\rightarrow Y$ the
natural map, and $i:D=\sigma^{-1}(Y)\rightarrow \tilde{X}_Y$ the
exceptional divisor of the blowing up. Then for each $p$, $k$ with
$k\geq 2p\geq 0$, we have the following isomorphism

$$
I_{p,k}: \bigg\{\bigoplus_{1\leq j \leq r-1}L_{p-j}H_{k-2j}(Y) \bigg\}\oplus
L_pH_k(X)\stackrel{\cong}{\longrightarrow} L_pH_k(\tilde{X}_Y)
$$
given by

$$
I_{p,k}(u_1,\cdots,u_{r-1},
u)=\sum_{j=1}^{r-1}i_*(h^j\cdot\pi^*u_j)+\sigma^*u.
$$
 }
\end{theorem}

\subsection{The Fulton-MacPherson configuration
spaces}\label{sec:2.2} {\hskip .2 in} Fulton and MacPherson have
constructed in \cite{FultonM} a compactification of the
configuration space of $n$ distinct labeled points in a non-singular
algebraic variety $X$. It is related to several areas of
mathematics. In their original paper, Fulton and MacPherson use it
to construct a differential graded algebra which is a model for
$F(X, n)$ in the sense of Sullivan ~\cite{FultonM}. Axelrod-Singer
constructed the compactification in the setting of smooth manifolds.
$\mathbb{P}^1[n]$ is related to the Deligne-Mumford compactification
$\overline{M}_{0,n}$ of the moduli space of nonsingular genus-$0$
projective curves.

Now we explain an explicit inductive construction of this
compactification given in \cite{FultonM}. $X[2]$ is the blowup of
$X^2$ along the diagonal $\Delta_{12}$. $X[3]$ is a sequence of
blowups of $X[2]\times X$ along non-singular subvarieties
corresponding to $\{\Delta_{123}; \Delta_{13}, \Delta_{23}\}$. More
specifically, denote by $\pi$ the blowup $X[2]\times X\to X^3$, we
blow up first along $\pi^{-1}(\Delta_{123})$, then along the strict
transforms of $\Delta_{13}$ and $\Delta_{23}$ (the two strict
transforms are disjoint, so they can be blown up in any order). In
general, $X[n+1]$ is a sequence of blowups of $X[n]\times X$ along
smooth subvarieties corresponding to all diagonals $\Delta_I$ where
$|I|\ge 2$ and $(n+1)\in I$.

Later, a symmetric construction of $X[n]$ has been given by several
people: De Concini and Procesi \cite{DP}, MacPherson and Procesi
\cite{MP}, and Thurston   \cite{Thurston}.  To construct $X[n]$ we
can blow up along diagonals by the order of ascending dimension,
which is different from the non-symmetric order of the original
construction. For example, $X[4]$ is the blowup of $X^4$ along
diagonals corresponding to:
$$1234;123,124,134,234;12,13,\dots,34.$$
Compare it with the order in \cite{FultonM}:
$$12;123;13,23;1234;124,134,234;14,24,34.$$

It is proved in \cite{author2} that, for any smooth projective
variety $X$ the following holds:

\begin{theorem}[\cite{author2}]
{ Let $X$ be a smooth projective variety defined
over $\mathbb{C}$. Then for each $p\ge 0$, there is an isomorphism
of Chow groups:
$${\rm Ch}^p(X[n])\cong\bigoplus_{\mathcal{S}}\bigoplus_{\underline{\mu}\in
M_{\mathcal{S}}} {\rm Ch}^{p-||\underline{\mu}||}(X^{c(\mathcal{S})}) .$$
where $\mathcal{S}$ runs through all nests of $[n]$.}
\end{theorem}

Notice that we use upper indices for the Chow groups in the above
theorem. By changing variable $\mu_I$ to $d(c_I-1)-\mu_I$, we get
exactly Theorem \ref{li thm} appeared in the introduction.

\begin{remark}
{The above theorem proved in \cite{author2} holds for non-singular
projective varieties $X$ over any algebraic closed field.}
\end{remark}

\medskip
Equivalently, but more explicitly, the Chow groups $X[n]$ of can be
calculated by using exponential generating functions. Here we adopt
R. Stanley's notation $[x^i]F(x)$ as the coefficient of $x^i$ in the
power series $F(x)$, which is generalized in an obvious way to the
following situation \cite{Stanley}:

\begin{equation}\label{eq:frac}
[\frac{x^it^n}{n!}]\sum_{j,q}a_{jq}\frac{x^jt^q}{q!}=a_{in}.
\end{equation}

\begin{corollary}\label{sec:cor2.1}
{If $h_i(x)$ are polynomials whose exponential generating
function $N(x,t)=\sum\limits_{i\ge 1}h_i(x)\frac{t^i}{i!}$ satisfies
the identity
$$(1-x)x^dt+(1-x^{d+1})=\exp(x^dN)-x^{d+1}\exp(N) ,$$
then  we have


$${\rm Ch}_p(X[n])=\bigoplus_{\substack{1\le m\le n\\
0\leq i\leq p}} {\rm Ch}_{p-i}(X^m)^{\oplus
[\frac{x^it^n}{n!}]\frac{N^m}{m!}}.$$

}
\end{corollary}

\qe

\subsection{Deligne-Beilinson cohomology} \hskip .2in
Let $X$ be a complex manifold of complex dimension $d$. Let
$\Omega_X^k$ the sheaf of holomorphic $k$-form on $X$. The
\textbf{Deligne complex of level p} is the complex of sheaves
$$\underline{\Z}_{\cal{D}}(p):0\rightarrow\Z\stackrel{(2i\pi)^p}{\rightarrow}\Omega_X^0\rightarrow
\Omega_X^1\rightarrow\Omega_X^2\rightarrow\cdots\rightarrow\Omega_X^{p-1}\rightarrow
0
$$
The \textbf{Deligne-Beilinson cohomology} of $X$ in level $p$ we
mean the hypercohomology of this complex:
$$ H^*_{\cal{D}}(X,\Z(p)):= \mathbb{H}^*(X,\underline{\Z}_{\cal{D}}(p))
$$

\medskip
There is a multiplication of complexes
$$ \nu: \Z(p)_{\cal{D}}\otimes \Z(q)_{\cal{D}}\rightarrow \Z(p+q)_{\cal{D}}
$$
defined as follows
$$
\nu(x\bullet y)=\left\{
\begin{array}{llll}
&x\cdot y, &&\quad {\rm if} \quad \deg ~x=0\\

&x\wedge dy, &&\quad {\rm if} \quad\deg ~x>0 \quad and\quad \deg ~y=q>0\\

&0, &&\quad {\rm otherwise}\\
\end{array}
\right.
$$
This gives a product structure on the Deligne--Beilinson cohomology
as follows

\begin{eqnarray}
\cup: H^k_{\cal{D}}(X,\Z(p))\otimes_{\Z}
H^{k'}_{\cal{D}}(X,\Z(q))\rightarrow H^{k+k'}_{\cal{D}}(X,\Z(p+q)).
\end{eqnarray}
For details, the reader is referred to \cite{Esnault-Viehweg}.

\medskip
Let $X$ be an $d$-dimensional compact K\"{a}hler manifold. The Hodge
filtration

$$\cdots\subseteq F^pH^k(X,{\mathbb{C}})\subseteq F^{p-1}H^k(X,{\mathbb{C}})
\subseteq\cdots\subseteq F^0H^k(X,{\mathbb{C}})=
H^k(X,{\mathbb{C}})$$ is defined by

$$ F^pH^k(X,{\mathbb{C}})=\oplus_{i\geq p}H^{i, k-i}(X).
$$

We denote by $p_X^k$ the natural quotient map  $p_X^k:
H^k(X,{\mathbb{C}})\rightarrow
H^k(X,{\mathbb{C}})/F^pH^k(X,{\mathbb{C}})$.

\medskip
It was proved (cf. \cite{Esnault-Viehweg}, Corollary 2.4;
\cite{Voisin1}, Proposition 12.26) that

{\footnotesize
\begin{eqnarray}
\cdots\rightarrow H^{k-1}(X,\mathbb{C})/F^pH^{k-1}(X,\mathbb{C})
\rightarrow H^k_{\cal{D}}(X,\Z(p)){\rightarrow} H^k(X,\Z)\rightarrow
H^k(X,\mathbb{C})/F^pH^k(X,\mathbb{C})\rightarrow \cdots
\end{eqnarray}
}

Now let $X$ be an $d$-dimensional projective variety over $\mathbb{C}$
and $i_0:Y\hookrightarrow X$ a smooth subvariety of codimension
$r\geq 2$. Let $\sigma:\tilde{X}_Y\rightarrow X$ be the blowup of
$X$ along $Y$, $\pi:D=\sigma^{-1}(Y)\rightarrow Y$ the natural map,
and $i:D=\sigma^{-1}(Y)\hookrightarrow \tilde{X}_Y$ the exceptional
divisor of the blowup. Set $U:= X-Y\cong \tilde{X}_Y - D$. Denote by
$j_0$ the inclusion $U\subset X$ and $j$ the inclusion $U\subset
\tilde{X}_Y$. Note that $\pi:D=\sigma^{-1}(Y)\rightarrow Y$ makes
$D$ into a projective bundle of rank $r-1$, given precisely by
$D=\P(N_{Y/X})$ and we have (cf. [\cite{Voisin2}, pg. 271])

$$ {\mathcal{O}}_{\tilde{X}_Y}(D)|_{D}=
{\mathcal{O}}_{\P(N_{Y/X})}(-1).
$$

Denote by $h$ the class of ${\mathcal{O}}_{\P(N_{Y/X})}(-1)$ under
the first Chern class $c_1: H^1(D,{\cal{O}}_D^* )\rightarrow H_{\cal
D}^2(D,\Z(1))$ (cf. [\cite{Esnault-Viehweg}, p. 88]).

\medskip
The following proposition was proved in \cite{Esnault-Viehweg}.

\begin{proposition}[\cite{Esnault-Viehweg}, Prop. 8.5]
{ The
Deligne-Beilinson cohomology $H^k_{\cal{D}}(D,\Z(p))$ of the
projective bundle $\pi: D\rightarrow Y$ is given by the following
isomorphism:
$$
\bigoplus_{0\leq j\leq r-1} H^{k-2j}_{\cal{D}}(Y,\Z(p-j))
\stackrel{\cong}{\longrightarrow}H^k_{\cal{D}}(D,\Z(p))
$$
}
\end{proposition}

\begin{remark}
{We omit the cup product of elements in
$H^{k-2j}_{\cal{D}}(Y,\Z(p-j))$ with $h^j$.}
\end{remark}

\medskip
Moreover, Barbieri-Viale proved the following blowup formula for
Deligne-Beilinson cohomology:

\begin{theorem}[\cite{Barbieri-Viale})]
 {\label{sec:BV} Let $X, Y, D , \tilde{X}_Y, Y$ be as above. Then for each
$p$, $k$ with $p\geq r\geq 0$, we have the following isomorphism

\begin{equation}
I_{p,k}: \bigg\{\bigoplus_{1\leq j \leq
r-1}H^{k-2j}_{\cal{D}}(Y,\Z(p-j)) \bigg\}\oplus
H^k_{\cal{D}}(X,\Z(p))\stackrel{\cong}{\longrightarrow}
H^k_{\cal{D}}(\tilde{X}_Y,\Z(p)).
\end{equation}
}
\end{theorem}

\begin{remark}
{Barbieri-Viale proved a general result, including the blowup
formula for \'{e}tale cohomology, to Theorem \ref{sec:BV}.}
\end{remark}


\section{Lawson homology for Fulton-MacPherson configuration spaces}

\hskip .2in In this section, we give a proof of Theorem
\ref{sec:M1}. According to the construction, the Fulton-MacPherson
configuration space $X[n]$ is obtained by a sequence of blowups
along all diagonals $\Delta_I$ in a suitable order. Each of
them is a blowup of a \emph{smooth} projective variety along a
\emph{smooth} projective subvariety. Therefore, we can calculate the
Lawson homology groups of $X[n]$ by successively applying the blowup formula
 for Lawson homology (Theorem \ref{sec:bulh}).

\medskip We have the following

\begin{theorem}\label{sec:3.1}
{\label{sec:LHFM} Let $X$ be a smooth projective variety defined
over $\mathbb{C}$. Then for each pair of integers $p,k$, $k\geq
2p\geq 0$, there is an isomorphism of Lawson homology groups:
$$
L_pH_k(X[n])\cong\bigoplus_{{\mathcal
S}}\bigoplus_{\underline{\mu}\in M_{\mathcal{S}}}
L_{p-||\underline{\mu}||}H_{k-2||\underline{\mu}||}(X^{c(\mathcal{S})}).
$$
}
\end{theorem}

\bp
This follows essentially from the construction of the
Fulton-MacPherson configuration space $X[n]$ and the blowup formula
for Lawson homology groups. The detailed computation for explicit
formulas will be given in the corollary below.

\qe

More explicitly, we have the following

\begin{corollary}\label{sec:cor3.1}
{Let $X$ be a smooth projective variety defined over
$\mathbb{C}$. Then for each pair of integers $p,k$, $k\geq 2p\geq
0$, there is an isomorphism of Lawson homology groups:
$$
L_pH_k(X[n])\cong\bigoplus_{\substack{1\le m\le n\\
0\leq i\leq p}} L_{p-i}H_{k-2i}(X^m)^{\oplus
[\frac{x^it^n}{n!}]\frac{N^m}{m!}}.
$$ where $N$  and $\oplus
[\frac{x^it^n}{n!}]\frac{N^m}{m!}$ are the same as those in
Corollary 2.1.
}
\end{corollary}

\medskip
\bp Let $h_n(x)$ be the polynomial
$$\sum_{\substack{\mathcal{S},
~\underline{\mu}\\c(\mathcal{S})=1}}x^{||\underline{\mu}||}.$$

Given a fixed nest $\mathcal{S}$ with $n$ leaves and
$c(\mathcal{S})=1$, its contribution to $h_n(x)$ is the product of
$\sigma_{c_I-1}$, where $I$ goes through all non-leaves of
$\mathcal{S}$ (if $\mathcal{S}$ has no non-leaves, i.e., it contains
only singletons, then the contribution is $1$). Therefore we have
the following recurrence formula
$$h_n(x)=\sum_{\{I_1,...,I_k\}\\\textrm{~partition of }
[n]}h_{|I_1|}h_{|I_2|}...h_{|I_k|}\sigma_{k-1}
.$$
where $\sigma_{k}=\sum_{i=1}^{dk-1}x^i$ for $k>0$, and $\sigma_0=0$.

By the Compositional Formula of exponential generating functions
(cf. \cite{Stanley}, Theorem 5.1.4), the generating function
$N(t):=\sum_{i\ge 1}h_i\frac{t^i}{i!}$ of $h_n$ satisfies the
identity
$$N-t+1=E_g(N),$$
where $E_g(N)=1+\sum_{i>0}\sigma_{i-1}N^i$.

Since $\sigma_j=(x^{jd}-x)/(x-1)$, calculation shows
$$E_g(N)=1+N+\frac{1}{x-1}\big{[}\frac{1}{x^d}(e^{x^dN-1}-1)-xe^N+x\big{]}.$$

Put it in the above identity, we have
$$(1-x)x^dt+(1-x^{d+1})=\exp(x^dN)-x^{d+1}\exp(N).$$

For a partition $\Pi=\{I_1,...,I_k\}$ of $[n]$, the number of times
of $h(\Delta_\Pi)(i)$ appear in the decomposition of $h(X[n])$ is
equal to $[x^k](h_{|I_1|}(x)...h_{|I_k|}(x))$. Add up this number
for all partitions with $k$ blocks, we will get the number of times
of $h(X^k)(i)$ appear in the decomposition, denoted by $a_{k,i}$.

Denote $$F_n(y)=\sum_{\{I_1,...,I_k\}\\\textrm{~partition of }
[n]}h_{|I_1}h_{|I_2}...h_{|I_k}y^k .$$ Then the coefficient
$[y^k]F_n(y)=\sum a_{k,i}x^i$. Use the Compositional Formula again,
$$F_n=[\frac{t^n}{n!}]\exp(yN).$$
Therefore
$$
\begin{array}{cc}
[y^k]F_n(y)&=[y^k][\frac{t^n}{n!}]\exp(yN)\\
&=[\frac{t^n}{n!}][y^k]\exp(yN)\\
&=[\frac{t^n}{n!}]\frac{N^k}{k!}.
\end{array}
$$

Now the results follow from Theorem \ref{sec:LHFM}.

\qe

\medskip
Similarly, we compute the  Deligne-Beilinson cohomology for Fulton-MacPherson
configuration spaces.

\begin{theorem}\label{sec:3.2}
{Let $X$ be a smooth projective variety defined over
$\mathbb{C}$. Then for each pair of integers $p$, $k$,  there is an
isomorphism of Deligne-Beilinson cohomology groups:
$$
H^k_{\cal{D}}(X[n],\Z(p))\cong\bigoplus_{\mathcal
S}\bigoplus_{\underline{\mu}\in M_{\mathcal{S}}}
H_{\cal{D}}^{k-2||\underline{\mu}||}\big(X^{c(\mathcal{S})},\Z({p-||\underline{\mu}||})\big).
$$
}
\end{theorem}

\medskip
\bp The method of the proof is the same as that in Theorem
\ref{sec:LHFM}. We get the result by using the explicit construction
of Fulton-MacPherson configuration spaces and Theorem
\ref{sec:BV}.\qe

\begin{remark}
{By using the same method, we can compute the \'{e}tale
cohomology for Fulton-MacPherson configuration spaces.}
\end{remark}

\begin{remark}
The decomposition of Lawson homology (Theorem \ref{sec:3.1}) and
Deligne - Beilinson cohomology (Theorem \ref{sec:3.2}) of the Fulton-MacPherson
configuration spaces can be generalized without any difficulty to
the wonderful compactifications of arrangements of subvarieties,
since the latter compactifications can also be constructed by a
sequence of blowups along smooth centers (for definition and
construction of these compactifications, see \cite{author2}).
\end{remark}

\section{Examples}

\begin{enumerate}
 \item
The Lawson homology group of $X[2]$.

The morphism $\pi: X[2]\to X^2$ is a blowup along the diagonal
$\Delta_{12}$. Corollary \ref{sec:cor3.1} asserts
$$L_pH_k(X[2])\cong
L_pH_k(X^2)\oplus\bigoplus_{j=1}^{d-1}L_{p-j}H_{k-2j}(X).$$

\item The Lawson homology group of $X[3]$.

Note that $X[3]$ is the blowup of $X^3$ first along small diagonal $\Delta_{123}$, then
along three disjoint proper transforms of diagonals $\Delta_{12}$,
$\Delta_{13}$ and $\Delta_{23}$.

Apply again Corollary \ref{sec:cor3.1}, we have
$$L_pH_k(X[3])\cong
L_pH_k(X^3)\oplus\bigoplus_{j=1}^{d-1}\big{(}L_{p-j}H_{k-2j}(X^2)\big{)}^{\oplus
3}$$
$$\oplus\bigoplus_{j=1}^{2d-1}\big{(}L_{p-j}H_{k-2j}(X)\big{)}^{\oplus
min\{3i-2, 6d-3i-2\}}
$$
\end{enumerate}

\section*{Acknowledge}

We would like to thank professor H. Blaine Lawson and Mark de Cataldo
for suggestions, conversations and all their helps.

\medskip

\noindent Wenchuan Hu, Department of Mathematics, MIT, Room 2-304,
77 Massachusetts Avenue Cambridge, MA 02139 \quad
Email: wenchuan@math.mit.edu

\medskip

\noindent Li Li, Department of Mathematics, Stony Brook University,
SUNY, Stony Brook, NY 11794-3651 \quad Email:lili@math.sunysb.edu


\begin{thebibliography}{AAAA}
%
%

\bibitem[Bv]{Barbieri-Viale}
L. Barbieri-Viale, {\em $ \mathcal{H}$-cohomologies versus algebraic
cycles.} (English. English summary) Math. Nachr. 184 (1997), 5--57.



\bibitem[DP]{DP}
C. De Concini and C. Procesi, {\em Wonderful models of subspace
arrangements}, Selecta Mathematica 1 (1995), 459--494.

\bibitem[EV]{Esnault-Viehweg} H. Esnault and E. Viehweg,
{\em Deligne-Be\u\i linson cohomology.} Be\u\i linson's conjectures
on special values of $L$-functions, 43--91, Perspect. Math., 4,
Academic Press, Boston, MA, 1988.


\bibitem[F]{Friedlander1} E. Friedlander, {\em Algebraic cycles, Chow
varieties, and Lawson homology.}  Compositio Math. 77 (1991), no. 1,
55--93.

\bibitem[FG]{FG} E. Friedlander and O. Gabber, {\em Cycle
spaces and intersection theory. Topological methods in modern
mathematics} (Stony Brook, NY, 1991), 325--370, Publish or Perish,
Houston, TX, 1993.

\bibitem[FrM]{Friedlander-Mazur} E. Friedlander and B. Mazur,
{\em Filtrations on the homology of algebraic varieties. With an
appendix by Daniel Quillen.}  Mem. Amer. Math. Soc. 110 (1994), no.
529, x+110 pp.

\bibitem[Fu]{Fulton}
W. Fulton, Intersection theory, Second edition, Springer-Verlag,
Berlin, 1998.

\bibitem[FuM]{FultonM}
W. Fulton and R. MacPherson, {\em A compactification of
configuration spaces}, Ann.\ Math.\ {\bf 139} (1994),
\hbox{183--225}.


\bibitem[GH]{Griffiths-Harris} Griffiths, Phillip; Harris, Joseph,  {\em Principles
of algebraic geometry.} Reprint of the 1978 original. Wiley Classics
Library.  John Wiley \& Sons, Inc., New York, 1994. xiv+813 pp. ISBN
0-471-05059-8


\bibitem[H]{author1} W. Hu, {\em Birational invariants defined by
Lawson homology.}\\  {arXiv:math.AG/0511722}.

\bibitem[L1]{Lawson1}
H. B. Lawson, {\em Algebraic cycles and homotopy theory.}, Ann. of
Math. {\bf 129}(1989), 253-291.

\bibitem[L2]{Lawson2} H. B. Lawson, {\em Spaces of algebraic
cycles.} pp. 137-213 in Surveys in Differential Geometry, 1995
vol.2, International Press, 1995.


\bibitem[Li]{author2} L. Li, {\em Chow motive of Fulton-MacPherson configuration
~spaces.} Ph. D. Thesis. SUNY, Stony Brook, 2006.

\bibitem[MP]{MP}
R. MacPherson and C. ~Procesi, {\em Making conical compactifications
wonderful}, Selecta Math. (N.S.) {\bf 4} (1998), no. 1, 125--139.


\bibitem[Pe]{Peters} C. Peters, {\em  Lawson homology for varieties with
small Chow groups and the induced filtration on the Griffiths
groups.}  Math. Z. 234 (2000), no. 2, 209--223.

\bibitem[St]{Stanley}
R. P. Stanley, {\em Enumerative combinatorics}, Vol. 2., Cambridge
Studies in Advanced Mathematics, 62. Cambridge University Press,
Cambridge, 1999.

\bibitem[Th]{Thurston}
D. Thurston, {\em Integral Expressions for the Vassiliev Knot
Invariants},\\ {arXiv:math.AG/9901110}.


\bibitem[V1]{Voisin1} C. Voisin, {\em Hodge theory and complex algebraic
geometry.} I. Translated from the French by Leila Schneps. Cambridge
Studies in Advanced Mathematics, 76. Cambridge University Press,
Cambridge, 2002. x+322 pp. ISBN 0-521-80260-1,

\bibitem[V2]{Voisin2} C. Voisin, {\em Hodge theory and complex algebraic
geometry.} II. Translated from the French by Leila Schneps.
Cambridge Studies in Advanced Mathematics, 77. Cambridge University
Press, Cambridge, 2003. x+351 pp. ISBN 0-521-80283-0


\end{thebibliography}
\end{document}